\documentclass[fleqn,10pt]{wlscirep}

\usepackage{graphicx,mathtools}
\usepackage{todonotes}
\RequirePackage{amsmath,amsfonts,amssymb,stmaryrd}
\RequirePackage{mathtools}
\RequirePackage{bm}
\RequirePackage{doi}
\RequirePackage{multirow}
\RequirePackage{siunitx}
\RequirePackage{seqsplit}

\graphicspath{{./figures/}}

\newcommand{\absperm}{\mathbb{K}}       % absolute permeability
\newcommand{\normal}{\bm{n}}        % normal vector
\newcommand{\porosity}{\phi}            % porosity
\newcommand{\head}{{h}}         % pressure head
\newcommand{\source}{{q}}           % source term
\renewcommand{\time}{{t}}           % time
\newcommand{\velocity}{\bm{v}}      % velocity
\newcommand{\vecu}{\bm{u}}

\newcommand{\dirichlet}{{\footnotesize \head}} % Dirichlet
\newcommand{\neumann}{{\footnotesize u}}       % Neumann
\newcommand{\K}{\absperm}

\renewcommand{\v}{\velocity}

\newcommand{\globalDomain}{\Omega}
\newcommand{\originalDomain}{\Lambda}

\newcommand{\pLn}{{\partial\originalDomain_{\neumann}}}
\newcommand{\pLd}{{\partial\originalDomain_{\dirichlet}}}

\newcommand{\jump}[1]{\left\llbracket {#1} \right\rrbracket}

\newcommand{\aperture}{\varepsilon}

\allowdisplaybreaks[1]
\title{Call for participation: Verification benchmarks for single-phase flow in three-dimensional fractured porous media}

\author[1,5]{Inga Berre}
\author[2]{Wietse Boon}
\author[2,*]{Bernd Flemisch}
\author[1]{Alessio Fumagalli}
\author[2]{Dennis Gläser}
\author[1]{Eirik Keilegavlen}
\author[3]{Anna Scotti}
\author[1]{Ivar Stefansson}
\author[4]{Alexandru Tatomir}
\affil[1]{Department of Mathematics, University of Bergen,
All\'{e}gaten 41, 5007 Bergen, Norway}
\affil[2]{Department of Hydromechanics and Modelling of Hydrosystems,
University of Stuttgart, Pfaffenwaldring 61, 70569 Stuttgart, Germany}
\affil[3]{Laboratory for Modeling and Scientific Computing MOX, Politecnico di Milano,
p.za Leonardo da Vinci 32, 20133 Milano, Italy}
\affil[4]{Department of Applied Geology, Geosciences Center,
University of G\"ottingen, Goldschmidtstrasse 3, 37077 G\"ottingen, Germany}
\affil[5]{Christian Michelsen Research, Bergen, Norway}

\affil[*]{bernd@iws.uni-stuttgart.de}

\begin{abstract}
This call for participation proposes four benchmark tests to verify and compare numerical schemes to solve single-phase flow in fractured porous 
media. With this, the two-dimensional suite of benchmark tests presented by Flemisch et al. \cite{Flemisch:2018:BSF} is extended to include three-dimensional problems. Moreover, transport simulations are included as a means to compare discretization methods for flow. 
With this publication, we invite researchers to contribute to the study by providing results to the test cases based on their applied discretization methods.
\end{abstract}
\begin{document}

\flushbottom
\maketitle
%thispagestyle{empty}
\section{Introduction}
In this call for participation, we propose four benchmark cases for the verification
of computational models for single-phase flow 
in three-dimensional, fractured porous media.
It can be seen as an extension of \cite{Flemisch:2018:BSF}, where
only two-dimensional test cases were considered. The current proposal also presents benchmarks where simulations of linear tracer transport are used as a means to compare flow discretizations.
However, the process of collecting and comparing results is different.
While the publication \cite{Flemisch:2018:BSF} contains both the
case descriptions and a comparison of several computational models, this time
we aim for a two-stage process. In particular, the present document
describes the benchmark cases and serves as a call for participation to the
interested scientific community. After a phase of collection and discussion of the results among the participants, we aim for a second publication
that presents a comparison of the participating models.
The case descriptions are accompanied by data in form of geometry descriptions,
simulation results and plotting scripts,
all available in the Git repository
\url{https://git.iws.uni-stuttgart.de/benchmarks/fracture-flow-3d.git}.
While access to this repository will be restricted to the benchmark participants
during the phase of collection and comparison of the results, all data will be made publicly
available upon submission of the comparison publication.

The rest of this document is organized as follows. In Section \ref{sec:model}, we describe
the mathematical models for fluid flow and transport in fractured porous media. Section \ref{sec:discretization} briefly describes the participating discretizations proposed by the authors as well as the discretization of the transport problem. Four test cases are introduced in Section
\ref{sec:cases} in form of domain and boundary condition descriptions as well as
pointers to the relevant data in the Git repository.
The final Section \ref{sec:participation} provides instructions on how to participate
in the benchmark study.

\section{Mathematical models}
\label{sec:model}

We introduce two models for flow and transport in fractured media. First, the flow model is presented in the conventional, equi-dimensional setting, allowing for a natural introduction to the physical parameters. From this formulation, we derive the mixed-dimensional model through appropriate reduction of the equations. Note that the mixed-dimensional model forms the focus of this study. Finally, we introduce the equi- and mixed-dimensional transport models.

\subsection{The equi-dimensional flow model}
\label{sec:equi_dimensional_flow_model}

We consider a steady-state incompressible single-phase flow through a porous medium
assumed to be described by Darcy's law, resulting in the governing system of equations
\begin{subequations}\label{eq:strong_dual}
\begin{gather}\label{eq:darcy}
	\begin{aligned}
 		\vecu + \K \nabla \head &= 0,  \\
 		\nabla \cdot \vecu &= \source,
    \end{aligned}
    \quad \text{in } \originalDomain.
\end{gather}
Here, $\vecu$ denotes the macroscopic
fluid velocity in \si{\metre\per\second} whereas $\K$ and $\head$ stand for hydraulic conductivity and
pressure head measured in \si{\metre\per\second} and \si{\metre}, respectively; $q$ represents a source/sink term measured in \si{\per\second}, and the regular domain $\originalDomain \subset \mathbb{R}^3$ will be called the
equi-dimensional domain. Coupled to \eqref{eq:darcy}, we consider boundary conditions on the boundary $\partial \originalDomain$ of
$\originalDomain$, namely
\begin{gather} \label{eq:darcybc}
	\begin{aligned}
 		\head |_{\pLd}&= \overline{\head} & \quad \text{on } \pLd, \\
		\vecu\cdot\normal |_{\pLn} &= \overline{u} & \text{on } \pLn.
	\end{aligned}
\end{gather}
\end{subequations}
We assume $\partial \originalDomain = \pLd \cup \pLn$, $\pLd \cap \pLn =
\emptyset$, and $|\pLd| > 0$. 
% WHY LAMNDA AND OMEGA FOR THE DOMAIN??
% ANSWER: we let \Lambda denote the equi-dimensional domain. \Omega then refers to the mixed-dimensional domain in which all fractures and intersections have been reduced in dimension. 
In \eqref{eq:darcybc} $\cdot |_A$ is a suitable trace operator on $A \subset \partial \originalDomain$, depending
on the quantity at hand. Also, $\overline{\head}$ indicates the pressure head imposed on the boundary $\pLd$, while
$\overline{u}$ is the prescribed Darcy velocity normal to the boundary $\pLn$ with respect to the outer unit normal vector $\normal$.

Problem \eqref{eq:strong_dual} can be recast in its primal formulation, obtaining the equations
\begin{gather}\label{eq:strong_primal}
	\begin{aligned}
	    -\nabla \cdot \K \nabla \head &=\source &\quad \text{in } \originalDomain, \\
    	\head |_{\pLd} &=  \overline{\head} & \text{on } \pLd, \\
		- \K \nabla \head \cdot\normal|_{\pLn} &=  \overline{u} & \text{on } \pLn.
    \end{aligned}
\end{gather}
Problem \eqref{eq:strong_dual} and \eqref{eq:strong_primal} are equivalent, however, different
numerical schemes are based on either of the two formulations. Under regularity assumptions on $\originalDomain$ 
and the data, 
the previous problems admit a  unique weak solution. We refer to \cite{Raviart1977,Brezzi1991,Roberts1991,Ern2004}
for more details. 

We assume that $\originalDomain$ contains several fractures, i.e., inclusions in the domain. The fracture walls are considered planar and smooth, and the fractures have two distinguishing features: (1) The fracture thickness, which we measure by the aperture, denoted by $\aperture$, is small compared to the extension of the fracture.
(2) The fracture hydraulic conductivity may differ significantly from that of the rest of $\originalDomain$, implying that the fractures may have significant impact on the flow in $\originalDomain$. 

We furthermore make the assumption that the principal directions of the local hydraulic conductivities  are aligned with the orientation of the fractures. In particular, the hydraulic conductivity in the matrix ($\K_3$), the fractures ($\K_2$), as well as in the intersections between two fractures ($\K_1$) and at the crossings of intersections ($\K_0$), can be decomposed in the following way:
\begin{align}
	\K_3 &= K_3^{eq}, &
    \K_2 &= \begin{bmatrix}
    	K_2^{eq} & \begin{matrix} 0 \\ 0 \end{matrix} \\
        \begin{matrix} 0 & 0 \end{matrix} & \kappa_2^{eq}
	\end{bmatrix}, & 
    \K_1 &= \begin{bmatrix}
		K_1^{eq} & 0 & 0 \\
        0 & \kappa_1^{eq} & 0 \\
        0 & 0 & \kappa_1^{eq} \\
	\end{bmatrix}, & 
    \K_0 &= \begin{bmatrix}
		\kappa_0^{eq} & 0 & 0 \\
        0 & \kappa_0^{eq} & 0 \\
        0 & 0 & \kappa_0^{eq} \\
	\end{bmatrix}.
\end{align} 
Here, $K_d^{eq}$ and $\kappa_d^{eq}$, for different values of $d$, denote the tangential and normal hydraulic conductivity, respectively. Thus, $K_d^{eq}$ is an elliptic $(d \times d)$-tensor function whereas $\kappa_d^{eq}$ is a positive scalar function. Note that the subscript $d$ indicates that the features will be represented by $d$-dimensional objects in the reduced model, as derived in the next section. The superscript $eq$, on the other hand, indicates that these quantities are related to the equi-dimensional model.

\subsection{Mixed-dimensional flow model}
\label{sec:mixed_dimensional_flow_model}
The small aperture  of the fractures justifies a reduction of dimensionality 
procedure to $\originalDomain$ where fractures and their intersections
are approximated by lower-dimensional objects. For more details on the derivation refer to \cite{Alboin2000,Faille2002,Angot2003,Martin2005,DAngelo2011,Fumagalli2012g,Schwenck2015,Flemisch2016,Boon2018,Fumagalli2017e}.

Here, we use $\globalDomain$ to denote the mixed-dimensional decomposition of $\originalDomain$. Let $\globalDomain$ with outer boundary $\partial \globalDomain$ be composed
of a possibly unconnected $3$-dimensional domain $\globalDomain_3$ which represents the rock matrix.
Furthermore, $\globalDomain$ contains up to $3$ lower-dimensional, open
subdomains, namely the fracture planes
$\globalDomain_{2}$, their intersection lines $\globalDomain_{1}$
and intersection points $\globalDomain_0$.
For compatibility, we assume that $\globalDomain_d \not\subset \globalDomain_{d'}$ for all $d' > d$. Finally, we introduce $\Gamma_{d} = \globalDomain_{d} \cap \partial \globalDomain_{d+1}$ as the set of $d$-interfaces
between inter-dimensional sub-domains $\globalDomain_{d}$ and $\globalDomain_{d+1}$ endowed with a normal unit vector $\normal$ pointing outward from $\globalDomain_{d+1}$. 

Remaining consistent with the notation convention above, data and unknowns will also be annotated with a subscript related to the dimension. As a first example, on a $d$-dimensional feature $\globalDomain_{d,i} \subseteq \globalDomain_d$ with counting index $i$, let $\aperture_{d,i}$ denote the cross-sectional volume, area, or length of the corresponding physical domain for $d = 0,..,2$ respectively. It naturally has the unit of measure \si{\metre\tothe{3-d}} and is extended as non-dimensional unity in $\globalDomain_3$. Moreover, we introduce for each $d$-feature with index $i$, a typical length $a_{d,i}$ such that $\aperture_{d,i} = a_{d,i}^{3 - d}$.
In the continuation, we will omit the subscript $i$ if no ambiguity arises.

We continue this subsection by first presenting the reduced model associated with \eqref{eq:strong_dual} in the two-dimensional fractures $\globalDomain_2$ followed by its generalization for all $d = 0, \dots, 3$.

\subsubsection{Two-dimensional fracture flow}

The primary variables in this formulation are the velocity $\vecu_3 = \vecu$ and hydraulic head $\head_3 = \head$ in the rock matrix $\globalDomain_3$, as well as the integrated, tangential velocity $\vecu_2$ and average hydraulic head $\head_2$ in the fracture. These are given pointwise for $x \in \globalDomain_2$ by
\begin{gather*}
	\vecu_2(x) = \int_{\aperture_2(x)} \bm{u}_{\|}
    \quad \text{and} \quad
    \head_2(x) = \dfrac{1}{\aperture_2(x)} \int_{\aperture_2(x)} \head.
\end{gather*}
Here, $\bm{u}_{\|}$ denotes the components of $\bm{u}$ tangential to $\Omega_2$. The integrals are computed in the normal direction of the fracture and thus the corresponding units of measurement are \si{\metre\squared\per\second} and \si{\metre} for $\vecu_2$ and $h_2$, respectively.

Let us derive the reduced Darcy's law by averaging and the mass balance equation by integration over the direction normal to the fractures. Recall that the vector $\normal$ here refers to the normal unit vector pointing from $\globalDomain_3$ into $\globalDomain_2$.
\begin{subequations}\label{eq:strong_fracture_dual 2}
  \begin{gather}\label{eq:darcy_fracture_dual 2}
      \begin{aligned}
          \frac{1}{\aperture_2} \vecu_{2} + K_2^{eq} \nabla_2 \head_2 &= 0  \\
          \nabla_2 \cdot \vecu_2 - \jump{\vecu_3 \cdot \normal} &= \source_2
      \end{aligned}
      \quad \text{in } \globalDomain_2,
  \end{gather}
with $\nabla_2$ the del-operator in the tangential directions and $\source_2$ representing the integrated source term, i.e. $\source_2(s) = \int_{\aperture_2(s)} \source$. Note that here, we have assumed $K_2^{eq}$ to be constant in the direction normal to $\globalDomain_2$. The jump operator is defined as
$\jump{\vecu_3 \cdot \normal} |_{\globalDomain_d} = \sum (\vecu_3 \cdot \normal |_{\Gamma_2})$, therewith representing the mass exchange between fracture and matrix. 
In particular, for each subdomain $\globalDomain_{2, i} \subseteq \globalDomain_2$, we sum over all flux contributions over sections of $\Gamma_{2}$ which coincide geometrically with $\globalDomain_{2, i}$.
These fluxes are assumed to satisfy the following Darcy-type law given by a finite difference between the hydraulic head in $\globalDomain_2$ and on $\Gamma_2$:
\begin{gather} \label{eq: Darcy normal 2}
	\vecu_3 \cdot \normal + \kappa_2^{eq} \frac{2}{a_d} (\head_2 - \head_3) = 0
    \quad \text{on } \Gamma_{2}.
\end{gather}
\end{subequations}
Note that to be mathematically precise, each term in this equation represents an appropriate trace or projection of the corresponding variable onto $\Gamma_2$.

\subsubsection{Generalized flow model}

Next, we generalize the equations described above to domains of all dimensions, thus including the intersection lines and points. For that purpose, we introduce the integrated velocity $\bm{u}_d$ for $d = 1$ and average hydraulic head $\head_d$ with $d = 0, 1$ given pointwise for $x \in \globalDomain_d$ by
\begin{gather*}
	\vecu_1(x) = \int_{\aperture_1(x)} \bm{u}_{\|}
    \quad \text{and} \quad
    \head_d(x) = \dfrac{1}{\aperture_d(x)} \int_{\aperture_d(x)} \head, \text{ for } d = 0, 1.
\end{gather*}
Again, $\bm{u}_{\|}$ denotes the components of $\bm{u}$ tangential to $\Omega_1$. The integrals are computed in the directions normal to the intersection lines $\Omega_1$, and points $\Omega_0$. The corresponding units of measurement are therefore \si{\metre\cubed\per\second}  and \si{\metre} for $\vecu_1$ and $h_d$, respectively. The analogues of \eqref{eq:darcy_fracture_dual 2} on these lower-dimensional manifolds are then given by
  \begin{gather}\label{eq:darcy_fracture_dual d}
      \begin{aligned}
          \frac{1}{\aperture_1} \vecu_{1} + K_1^{eq} \nabla_1 \head_1 &= 0  \\
          \nabla_1 \cdot \vecu_1 - \jump{\vecu_2 \cdot \normal} &= \source_1
      \quad \text{in } \globalDomain_1, \\
      - \jump{\vecu_1 \cdot \normal} &= \source_0
      \quad \text{in } \globalDomain_0.
      \end{aligned}
  \end{gather}
Here, $\nabla_1$ denotes the del-operator, i.e. the derivative, in $\Omega_1$. Moreover, the linear jump operator $\jump{\cdot}$ is naturally generalized to
$\jump{\vecu_{d+1} \cdot \normal} |_{\globalDomain_d} = \sum (\vecu_{d+1} \cdot \normal |_{\Gamma_{d}})$, where we for each subdomain $\globalDomain_{d, i} \subseteq \globalDomain_d$ sum over all flux contributions over sections of $\Gamma_{d}$ which coincide geometrically with $\globalDomain_{d, i}$.  Finally $q_1$ and $q_0$ correspond to the integrated source terms in the intersection lines and points, respectively.

Due to our choice of defining $\bm{u}_d$ as the integrated velocity, a scaling with $\aperture_{d + 1}$ appears in the equation governing the flux across $\Gamma_d$:
\begin{gather} \label{eq: Darcy normal d}
	\frac{1}{\epsilon_{d + 1}} \vecu_{d + 1} \cdot \normal + \kappa_d^{eq} \frac{2}{a_d} (\head_d - \head_{d + 1}) = 0
    \quad \text{on } \Gamma_{d}, \ d = 0, 1.
\end{gather}
Recalling that $\epsilon_3 = 1$, it now follows that the effective tangential and normal hydraulic conductivities are given by:
\begin{subequations}
	\begin{align}
      K_d &= \aperture_d K_d^{eq}, 
      &\text{in } \globalDomain_d, \ d &= 1, \dots, 3 \\
      \kappa_d &=  \aperture_{d + 1} \frac{2}{a_d} \kappa_d^{eq}, 
      &\text{on } \Gamma_d, \ d &= 0, \dots, 2.
	\end{align}
\end{subequations}
From these definitions, it is clear that the units of $K_d$ and $\kappa_d$ are \si{\metre\tothe{4-d}\per\second} and \si{\metre\tothe{2-d}\per\second}, respectively.

Collecting the above equations, we obtain the generalization of system \eqref{eq:strong_fracture_dual 2} to subdomains of all dimensions. The resulting system consists of Darcy's law in both tangential and normal directions followed by the mass conservation equations: 
\begin{subequations}\label{eq:strong_fracture_dual detailed}
	\begin{align}
 		\vecu_d + K_d \nabla_d \head_d &= 0, 
        &\text{in } &\globalDomain_d, \ d = 1, \dots, 3,\\
		\vecu_{d+1} \cdot \normal + \kappa_d (\head_d - \head_{d+1}) &= 0,
    	&\text{on } &\Gamma_d, \ d = 0, \dots, 2,\\
        \nabla_d \cdot \vecu_3 &= \source_3,
    	&\text{in } &\globalDomain_3, \\
 		\nabla_d \cdot \vecu_d - \jump{\vecu_{d+1} \cdot \normal} &= \source_d,
    	&\text{in } &\globalDomain_d, \ d = 1, 2,\\
 		- \jump{\vecu_1 \cdot \normal} &= \source_0,
    	&\text{in } &\globalDomain_0.
    \end{align}
\end{subequations}
The source term is given by $\source_3$ for the rock matrix and $\source_d(x) = \int_{\aperture_d(x)} \source$ measured in \si{\metre\tothe{3-d}\per\second}.

The system \eqref{eq:strong_fracture_dual detailed} is then compactly described by:
\begin{subequations}\label{eq:strong_fracture_dual}
	\begin{align}
 		\vecu_d + K_d \nabla_d \head_d &= 0, 
        &\text{in } \globalDomain_d, \ d &= 1, \dots, 3, \label{eq:strong_Darcy_t} \\
		\vecu_{d+1} \cdot \normal + \kappa_d (\head_d - \head_{d+1}) &= 0,
    	&\text{on } \Gamma_d, \ d &= 0, \dots, 2, \label{eq:strong_Darcy_n} \\
 		\nabla_d \cdot \vecu_d - \jump{\vecu_{d+1} \cdot \normal} &= \source_d,
    	&\text{in } \globalDomain_d, \ d &= 0, \dots, 3. \label{eq:strong_massconv}
    \end{align}
\end{subequations}
in which the non-physical $\vecu_4$ and $\vecu_0$ are understood as zero. 
%\as{ what is the meaning of u4?} FIXED
The associated boundary conditions are inherited from the equidimensional model with the addition of a no-flux condition at embedded fracture endings:
\begin{subequations}\label{eq:strong_fracture_dual_BC}
	\begin{align}
        \head_d &= \overline{\head} 
        & \text{on } \partial \globalDomain_d \cap \pLd, \ d &= 0, \dots, 3,\\
		\vecu_d \cdot \normal &= \aperture_d \overline{u} 
        & \text{on } \partial \globalDomain_d \cap \pLn, \ d &= 1, \dots, 3, \\
        \vecu_d \cdot \normal &= 0 
        & \text{on } \partial \globalDomain_d \backslash (\Gamma_{d - 1} \cup \partial \originalDomain), \ d &= 1, \dots, 3.
    \end{align}
\end{subequations}

To finish the section, we present the primal formulation of the mixed-dimensional fracture flow model. In analogy with \eqref{eq:strong_primal}, this formulation is derived by substituting the Darcy's laws \eqref{eq:strong_Darcy_t}, \eqref{eq:strong_Darcy_n} into the conservation equation \eqref{eq:strong_massconv}:
\begin{align}\label{eq:strong_fracture_primal}
	-\nabla_d \cdot K_d \nabla_d \head_d + \jump{\kappa_d (\head_d - \head_{d+1})} 
    &= \source_d,
   	&\text{in } \globalDomain_d, \ d &= 0, \dots, 3.
\end{align}
In this case, we set the divergence term to zero if $d = 0$ and the jump term to zero if $d = 3$. The boundary conditions are given by
\begin{subequations}\label{eq:strong_fracture_primal_BC}
\begin{align}
        \head_d &= \overline{\head} 
        & \text{on } \partial \globalDomain_d \cap \pLd, \ d &= 0, \dots, 3,\\
		- K_d \nabla_d \head_d \cdot \normal &= \aperture_d \overline{u} 
        & \text{on } \partial \globalDomain_d \cap \pLn, \ d &= 1, \dots, 3, \\
        - K_d \nabla_d \head_d \cdot \normal &= 0 
        & \text{on } \partial \globalDomain_d \backslash (\Gamma_{d - 1} \cup \partial \originalDomain), \ d &= 1, \dots, 3.
\end{align}
\end{subequations}

Many discretization schemes of the mixed-dimensional model described here, ignore flow in one-dimensional fracture intersections and zero-dimensional intersections thereof. Although these correspond to discretizing a simpler model, this is perfectly in line with the proposed study.
In particular, it will be interesting to evaluate differences in the results that can be attributed to whether these lower-dimensional intersections are included in the model, and therewith the discretization, or not.

\subsection{Equi-dimensional transport model}
\label{sec:equi_dimensional_transport_model}
We now consider a scalar quantity $c$ with the unit of measure \si{\metre\tothe{-3}}, which is transported through the porous medium subject to the velocity field resulting from the flow model presented in the previous sections. The purely advective transport of $c$ is described by the conservation equation:
\begin{equation}\label{eq:strong_transport}
	\porosity \frac{\partial c}{\partial \time} + \nabla \cdot \left( c \vecu \right) = \source_c	\quad \text{in } \originalDomain,
\end{equation}
where $\porosity$ is the porosity of the porous medium and $\source_c$ is a source/sink term for $c$ given in \si{\metre\tothe{-3}\per\second}. We define Dirichlet boundary conditions on those boundary segments where inflow occurs, i.e.:
\begin{equation} \label{eq:transportbc}
 		c |_{\partial\originalDomain_{\footnotesize c}}
            = \overline{c} \quad \text{on } \partial\originalDomain_{\footnotesize c}, \quad
              \partial\originalDomain_{\footnotesize c} =
                 \{ x \in \partial\originalDomain: \,\, \vecu \cdot \normal < 0 \},
\end{equation}
with $\overline{c}$ being the value for $c$ prescribed on the boundary $\partial\originalDomain_{\footnotesize c}$.

\subsection{Mixed-dimensional transport model}
\label{sec:mixed_dimensional_transport_model}
In analogy to Section \ref{sec:mixed_dimensional_flow_model}, we choose the average value for $c$ as primary variable, which is defined as $c_3 = c$ in $\globalDomain_3$ and for the lower dimensional objects (with $d \le 2$) as
\begin{equation*}
c_d(s) = \dfrac{1}{\aperture_d(s)} \int_{\aperture_d(s)} c.
\end{equation*}
Following the derivation of the mixed-dimensional flow model presented in Section \ref{sec:mixed_dimensional_flow_model}, the resulting mixed-dimensional transport model reads:
\begin{equation}\label{eq:strong_transport_dual}
\aperture_d \porosity_d \frac{\partial c_d}{\partial \time}
		 + \nabla_d \cdot \left( c_d \vecu_d \right)
		 - \jump{\tilde{c}_{d+1} \left( \vecu_{d+1} \cdot \normal \right)}
		 = \source_{c,d}	\quad \text{in } \globalDomain_d, 
           \quad d = 0, \dots, 3.
\end{equation}
Note that for $d=0$, the divergence term is void. 
Here, the porosity is simply $\porosity_d = \porosity^{eq}$, with unit of measure \si{\metre\tothe{-3}}, and $\tilde{c}_{d+1}$ is evaluated on the basis of a first order upwind scheme, i.e.:
\begin{equation} \label{eq:upwind}
  \tilde{c}_{d+1} =
        \begin{cases}
            c_{d+1} & \text{if} \quad \vecu_{d+1} \cdot \normal |_{\Gamma_{d}} > 0\\
            c_d   & \text{if} \quad \vecu_{d+1} \cdot \normal |_{\Gamma_{d}} < 0.
        \end{cases}
\end{equation}
As in the flow model, the jump operator represents the sum of the fluxes over all contributions defined on sections of $\Gamma_{d}$ that coincide geometrically with $\globalDomain_{d, i}$.

\section{Discretization} \label{sec:discretization}
% WMB: Used to be a subsection, but it doesn't belong in "mathematical models" IMHO.
The intent of the proposed benchmark study is to quantitatively evaluate different
discretization schemes for the mixed-dimensional flow models \eqref{eq:strong_fracture_dual},
\eqref{eq:strong_fracture_dual_BC} and \eqref{eq:strong_fracture_primal},
\eqref{eq:strong_fracture_primal_BC}.
The transport models described above are employed as a means to further evaluate and compare results for different flow discretizations. In the repository, the authors have already uploaded results to the test cases based on the following methods:
% \as{(The invited researchers needs to get to know which methods that are already covered, so I think we need to put this information somewhere)} 
% WMB: Agreed, these descriptions need to include the defining features. I believe all methods have a mixed-dimensional representation of fractures?
\begin{enumerate}
\item \texttt{UiB-TPFA}: Two-point flux approximation finite volume method (cell-centered pressure, conforming DFM grid, discontinuous over lower-dimensional features, mass transfer in intersection lines and points).
\item \texttt{UiB-MPFA}: Multi-point flux approximation finite volume method (cell-centered pressure, conforming DFM grid, discontinuous over lower-dimensional features, mass transfer in intersection lines and points).
\item \texttt{UiB-RT0}: Lowest order Raviart-Thomas mixed finite elements (cell-centered pressure and face flux, conforming DFM grid, discontinuous over lower-dimensional features, mass transfer in intersection lines and points).
\item \texttt{UiB-MVEM}: Lowest order mixed virtual element method (cell-centered pressure and face flux, conforming DFM grid, discontinuous over lower-dimensional features, mass transfer in intersection lines and points).
\item \texttt{USTUTT-MPFA}: Multi-point flux approximation finite volume method (cell-centered pressure, conforming DFM grid, discontinuous over lower-dimensional features, \emph{no} mass transfer in intersection lines and points).
%\item Vertex-centred, continuous-pressure, conforming lower-dimensionalDFM (Box)
%\item Continuous-pressure, non-conforming embedded DFM (EDFM)
%\item Discontinuous-pressure, non-conforming dual XFEM (D-XFEM)
\end{enumerate}
Methods 1-4 are all implemented following the mixed-dimensional approach described in \cite{Nordbotten2018}. 

To evaluate the fluxes produced by each discretization quantitatively, they are to be inserted into a standard
cell-centered first-order-upwind scheme for the transport, together with the implicit Euler method for temporal discretization with a fixed time-step that will be prescribed for each test case. 

\section{Benchmark cases}
\label{sec:cases}
In the next subsections we introduce the proposed benchmark cases. For each case, 
pressure and tracer concentration need to be computed along with several associated
macroscopic metrics. In Subsection \ref{subsec:single_fracture} a problem containing a single fracture problem is
considered. Subsection \ref{subsec:regular} contains a synthetic 
network composed by 9 fractures. The example in Subsection \ref{subsec:small_features} considers
the geometrically challenging case of almost intersecting fractures, fractures with small intersections, and
other features that a fracture network may exhibit. Finally, in Subsection
\ref{subsec:field_network} we study a selection of a few dozen fractures from a real network.
In all the cases the purpose is to validate the performance of the proposed numerical schemes and model choices.

\newpage 
\subsection{Case 1: single fracture}\label{subsec:single_fracture}

\begin{figure}[!b]
	\centering
    \includegraphics[width=0.5\textwidth]{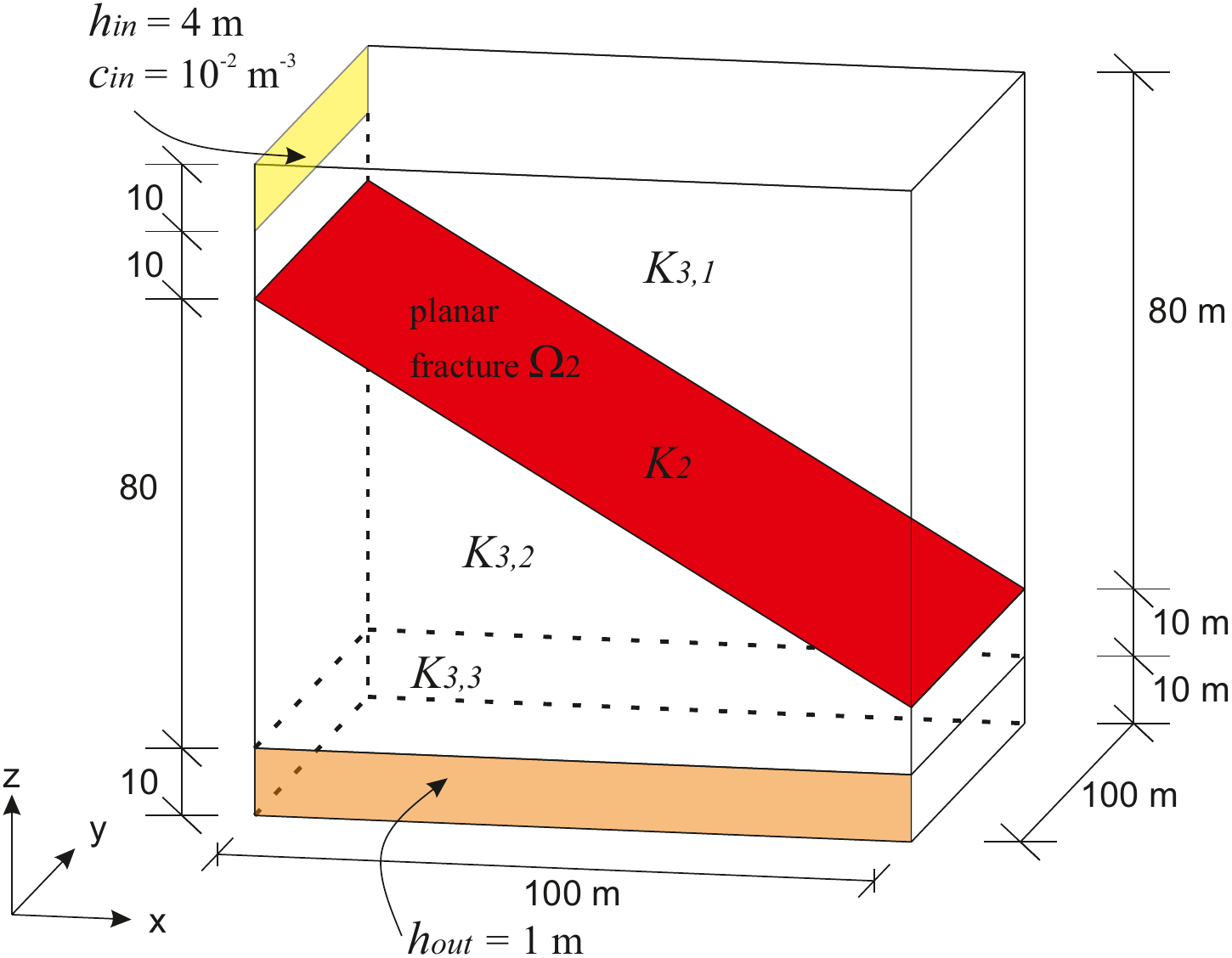}
    \caption{Conceptual model and geometrical description of the domain for test case of Subsection \ref{subsec:single_fracture}.}
    \label{fig:problem1}
\end{figure}

\paragraph{Description}
Figure \ref{fig:problem1} illustrates the first proposed benchmark example for testing the discrete fracture flow and transport models in three-dimensional space. 
The geometry of the system is slightly modified from the previous works of \cite{Zielke:1991:DMT} and \cite{Barlag:1998:AMM}. 
The domain, $\Omega$, is a cube-shaped region (\SI{0}{\metre}, \SI{100}{\metre}) x (\SI{0}{\metre}, \SI{100}{\metre}) x (\SI{0}{\metre}, \SI{100}{\metre}) which is crossed by a planar fracture, $\globalDomain_2$, with a thickness of \SI{1e-2}{\metre}. The fracture divides the domain into the sub-domains $\globalDomain_{3,1}$ and $\globalDomain_{3,2}$. Furthermore, an additional sub-domain $\globalDomain_{3,3}$ is defined in the lower part of the domain, representing a heterogeneity within the rock matrix. Inflow into the system occurs through a narrow band defined by $\{\SI{0}{\metre}\}\times (\SI{0}{\metre}, \SI{100}{\metre}) \times (\SI{90}{\metre}, \SI{100}{\metre}) $. Similarly, the outlet is a narrow band defined by $(\SI{0}{\metre}, \SI{100}{\metre}) \times \{\SI{0}{\metre}\}\times (\SI{0}{\metre}, \SI{100}{\metre})$. Except for the inlet and outlet bands, on which Dirichlet boundary conditions are set, all other boundaries are treated as no-flow boundaries. For the hydraulic head, $h_{in} = \SI{4}{\metre}$ and $h_{out} = \SI{1}{\metre}$ are set on the inlet and outlet boundaries, respectively, while for the transport $c_{in} = \SI{1e-2}{\metre\tothe{-3}}$ is set at the inlet. The overall simulation time is \SI{1e9}{\second}.

\paragraph{Parameters}
\begin{center}
\begin{tabular}{|l|ll|}\hline
%\textbf{Parameter}                         & \textbf{Value} &         \\\hline
Matrix hydraulic conductivity $K_{3,1}$, $K_{3,2}$ & \num{1e-6}$\bm{I}$ & \si{\metre\per\second} \\
Matrix hydraulic conductivity $K_{3,3}$         & \num{1e-5}$\bm{I}$ & \si{\metre\per\second} \\
Fracture effective tangential hydraulic conductivity $K_2$ & \num{1e-3}$\bm{I}$ & \si{\metre^2\per\second} \\
Fracture effective normal hydraulic conductivity $\kappa_2$ & \num{20} &\si{\per\second} \\
Matrix porosity $\phi_{3,1}, \phi_{3,2}$         & \num{2e-1}                              &\\
Matrix porosity $\phi_{3,3}$                   & \num{2.5e-1}                             & \\
Fracture porosity $\phi_2$                   & \num{4e-1}                              & \\
Fracture cross-sectional length $\epsilon_2$                        & \num{1e-2}                      & \si{\metre} \\\hline
\end{tabular}
\end{center}
%\caption{Parameter values for test case of Subsection \ref{subsec:single_fracture}.}
%\label{tb:tab1}

\paragraph{Folder in Git repository} \texttt{single}

\paragraph{Results} 
We consider a constant time step size $\Delta t = \SI{e7}{\second}$. The results of three different simulations with approximately $1k$, $10k$ and $100k$ cells for the 3D domain are collected
and the comparison among the different schemes is then done on the basis of the following results:
\begin{enumerate}
    \item The number of cells and degrees of freedom, according to the guidelines described in Subsection \ref{subseq:data_reporting} point \ref{data:results}, saved in the file \texttt{results.csv};
    
    \item \label{case1_result1} the computation of  $\int_{\globalDomain_{3,3}} \porosity_3 c_3 \, \mathrm{d}x$  for each time step;
    
    \item \label{case1_result2} the computation of $\int_{\globalDomain_2} \aperture_2 \porosity_2 c_2 \, \mathrm{d}x$  for each time step;
    
    \item \label{case1_result3} the integrated flux of $c$ across the outlet boundary for each time step;
    
    \item \label{case1_result4} the plot of the hydraulic head $\head_3$ in the matrix along the line $\left( \SI{0}{\metre}, \SI{100}{\metre}, \SI{100}{\metre} \right)$-$\left( \SI{100}{\metre}, \SI{0}{\metre}, \SI{0}{\metre} \right)$;
    
    \item \label{case1_result5} the plot of $c_3$ in the matrix, at the final simulation time, along the line $\left( \SI{0}{\metre}, \SI{100}{\metre}, \SI{100}{\metre} \right)$-$\left( \SI{100}{\metre}, \SI{0}{\metre}, \SI{0}{\metre} \right)$;
    
    \item \label{case1_result6} the plot of $c_2$ within the fracture, at the final simulation time, along the line $\left( \SI{0}{\metre}, \SI{100}{\metre}, \SI{80}{\metre} \right)$-$\left( \SI{100}{\metre}, \SI{0}{\metre}, \SI{20}{\metre} \right)$.
    
\end{enumerate}
In files called \texttt{dot\_refinement\_\$\{REFINEMENT\}.csv}
report the results ordered column-wise in the following way: time, result from point \ref{case1_result1}, \ref{case1_result2}, and \ref{case1_result3}. The variable \texttt{\$\{REFINEMENT\}} runs from 0 to 2 for increasing refinement. Similarly, in files called \texttt{dol\_refinement\_\$\{REFINEMENT\}.csv} report the results ordered in the following way: arc length and result from point \ref{case1_result4}, arc length and result from point \ref{case1_result5}, and arc length and result from point \ref{case1_result6}.
See Subsection \ref{subseq:data_reporting} for more details on how to report results.

\newpage

%%%%%%%%%%%%%%%%%%%%%%%%%%%%%%%%%%%%%%%%%%%%%%%%%%%%%%%%%%%%%%%%%%%%%%%%%%%%%%%%%%%%%%%%%%%%%%%%%%%%%%%%%%%%%%%%%%%%%%%%%%

\subsection{Case 2: regular fracture network}\label{subsec:regular}

\begin{figure}[!b]
	\centering
    \includegraphics[width=0.33\textwidth]{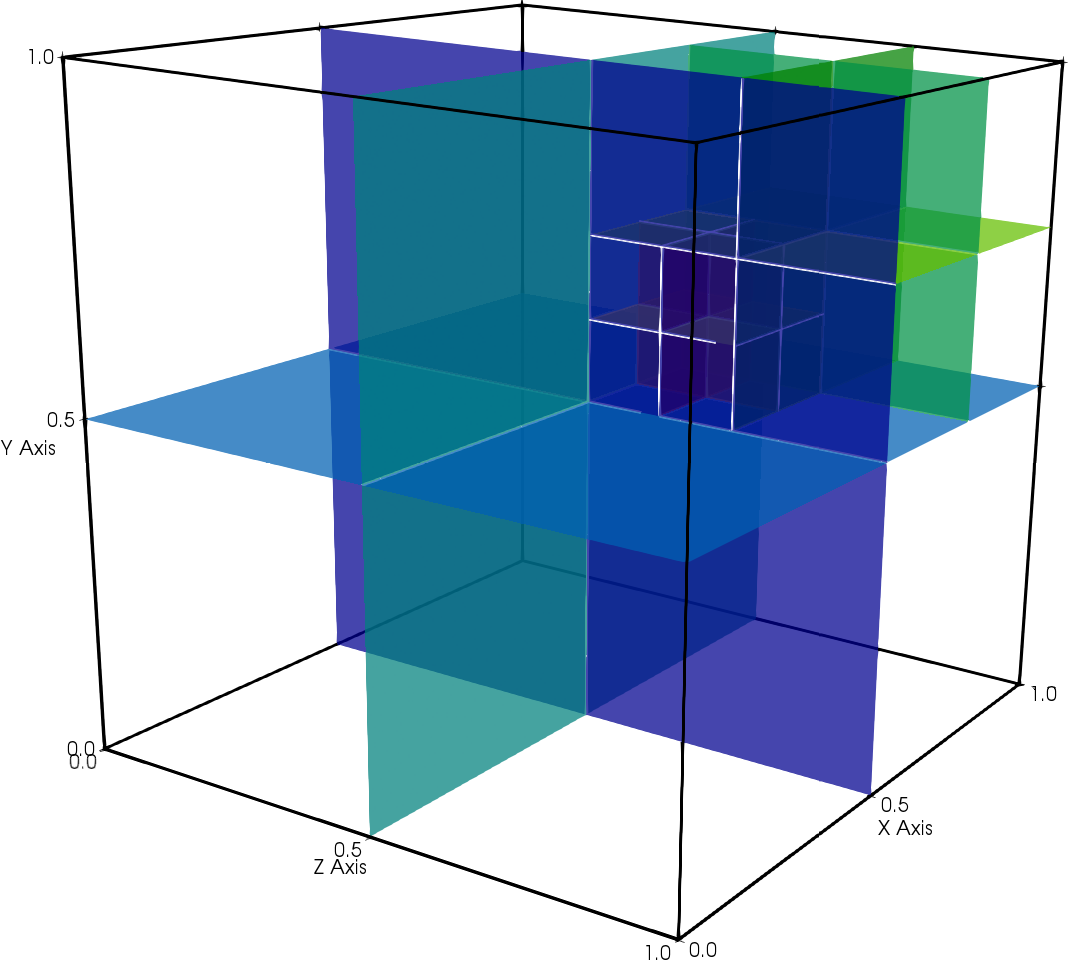}%
    \hspace*{0.1\textwidth}%
    \includegraphics[width=0.33\textwidth]{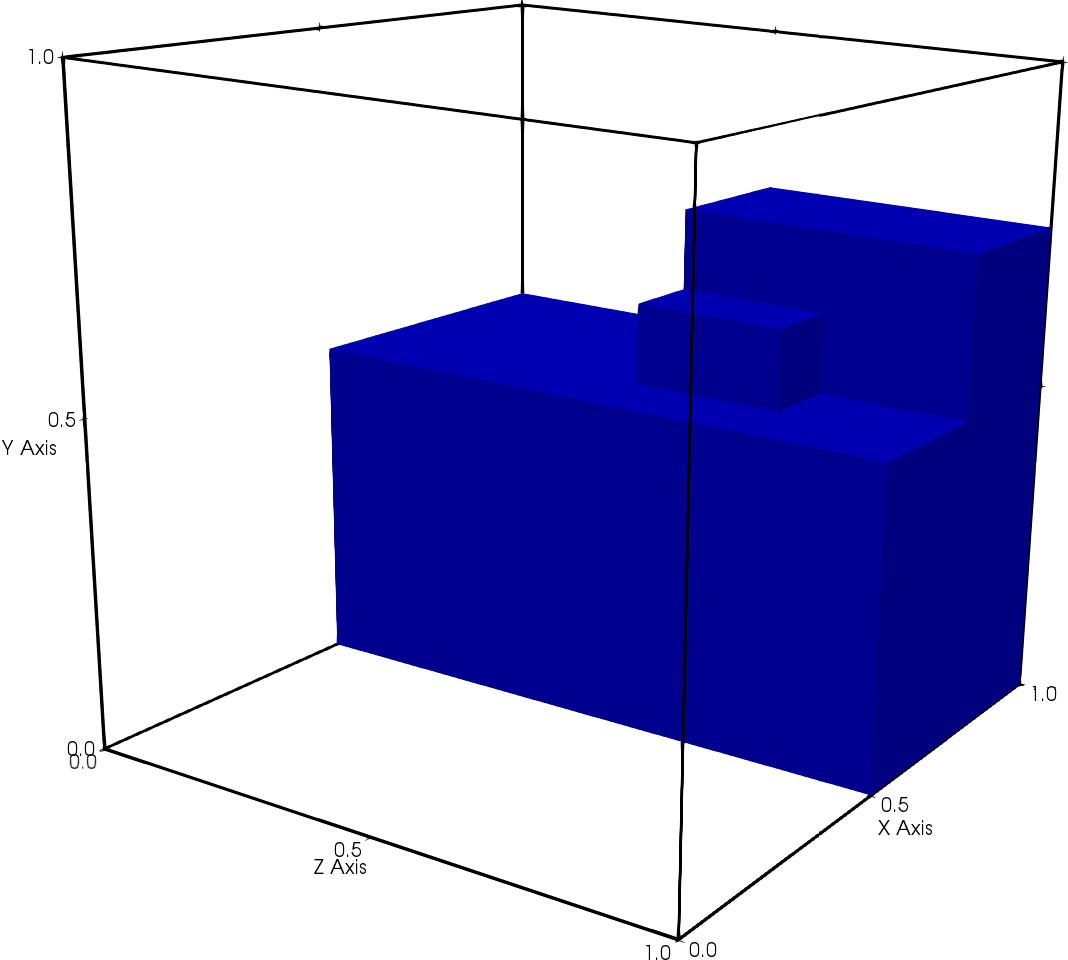}
    \caption{On the left, representation of the fractures and the outline of the domain
    	for the test case of Subsection \ref{subsec:regular}. On the right, the part of the 
        rock matrix with low hydraulic conductivity $\globalDomain_{3,1} $.}
    \label{fig:domain_regular}
\end{figure}

\paragraph{Description}
The second benchmark is an extension of test case 4.1 from the 
benchmark study \cite{Flemisch:2018:BSF} in a 3D setting. The geometry is slightly modified 
from \cite{Geiger:2011:NMD} to allow a full 
3D solution field. The unit cube domain is $\Omega = \left(\SI{0}{\metre}, \SI{1}{\metre} \right)^3$ and 9 fractures are present, 
see Appendix \ref{appendix:regular} for the coordinates. A graphical representation is given 
by Figure \ref{fig:domain_regular}.

We define the portion of the boundary where we impose a
non-homogeneous flux condition as $\partial \globalDomain_{flux}$ and choose
$\partial \globalDomain_{h} = \{ (x, y, z) \in \partial \Omega: x, y, z > \SI{0.875}{\metre}\}$ and
$\partial \globalDomain_{flux} = \{ (x, y, z) \in \partial \Omega: x, y, z < \SI{0.25}{\metre}\}$, respectively.
On $\partial \globalDomain_{h}$ we impose $\overline{h} = \SI{1}{\metre}$, while on $\partial \globalDomain_{flux}$ we impose $\overline{u} = \SI{-1}{\metre\per\second}$.
On the rest of the boundary of $\Omega$ and on the lower-dimensional objects, we impose homogeneous flux conditions. Two test cases are considered: high-conductive fractures (\texttt{\$\{COND\} = 0}) and 
low-conductive fractures  (\texttt{\$\{COND\} = 1}). We define the following matrix sub-regions 
\begin{align*}
    \globalDomain_{3,0} =&\ \globalDomain_{3} \setminus \globalDomain_{3,1} \\
    \globalDomain_{3,1} =&\ \{(x, y, z) \in \globalDomain_{3}: x > \SI{0.5}{\metre} \cap y < \SI{0.5}{\metre}\} \\
    &\cup \{(x, y, z) \in \globalDomain_{3}: x > \SI{0.75}{\metre} \cap \SI{0.5}{\metre} < y < \SI{0.75}{\metre} \cap z > \SI{0.5}{\metre}\}\\
    &\cup \{(x, y, z) \in \globalDomain_{3}: \SI{0.625}{\metre} < x < \SI{0.75}{\metre} \cap \SI{0.5}{\metre} < y < \SI{0.625}{\metre} \cap \SI{0.5}{\metre} < z < \SI{0.75}{\metre}\},
\end{align*}
see the right part of Figure \ref{fig:domain_regular}, and prescribe two different
matrix hydraulic conductivities.

For the transport problem, we impose unitary 
concentration at the inflow boundary and time interval equal to $[0, 0.25]$ \si{\second}.

\paragraph{Parameters}
\begin{center}
\begin{tabular}{|l|ll|ll|}\hline
                        & \multicolumn{2}{l|}{\texttt{\$\{COND\} = 0}} &  \multicolumn{2}{l|}{\texttt{\$\{COND\} = 1}}      \\\hline
Matrix hydraulic conductivity $K_{3}|_{\Omega_{3, 0}}$        & $\bm{I}$& \si{\meter\per\second} & $\bm{I}$&\si{\meter\per\second}\\
Matrix hydraulic conductivity $K_{3}|_{\Omega_{3, 1}}$         & \num{1e-1}$\bm{I}$&\si{\meter\per\second}& \num{1e-1}$\bm{I}$&\si{\meter\per\second} \\
Fracture effective tangential hydraulic conductivity $K_2$ & $\bm{I}$ & \si{\meter\squared\per\second} & \num{1e-8}$\bm{I}$ & \si{\meter\squared\per\second}\\
Fracture effective normal hydraulic conductivity $\kappa_2$ & \num{2e8} & \si{\per\second} & $2$& \si{\per\second}\\
Intersection effective tangential hydraulic conductivity $K_1$ & \num{1e-4} &\si{\meter\cubed\per\second} & \num{1e-12} & \si{\meter\cubed\per\second}\\
Intersection effective normal hydraulic conductivity $\kappa_1$ & \num{2e4} &\si{\meter\per\second} & \num{2e-4} & \si{\meter\per\second}\\
Intersection effective normal hydraulic conductivity $\kappa_0$ & \num{2} &\si{\meter\squared\per\second} & \num{2e-8}& \si{\meter\squared\per\second}\\
Matrix porosity $\phi_{3}$ & \num{1e-1} & & \num{1e-1} & \\
Fracture porosity $\phi_2$ & \num{9e-1} & & \num{1e-2}& \\
Intersection porosity $\phi_1$ & \num{9e-1} & & \num{1e-2}& \\
Fracture cross-sectional length $\epsilon_2$  & \num{1e-4} &\si{\meter}& \num{1e-4} &\si{\meter}\\
Intersection cross-sectional area $\epsilon_1$  & \num{1e-8} &\si{\meter\squared} & \num{1e-8} &\si{\meter\squared}\\
Intersection cross-sectional volume $\epsilon_0$ & \num{1e-12} & \si{\meter\cubed} & \num{1e-12} &\si{\meter\cubed}\\
\hline
\end{tabular} 
\end{center}

\paragraph{Folder in Git repository} \texttt{regular}

\paragraph{Results}
Results will be collected for a sequence of 3 simulations, by taking approximately $500$, $4k$, and $32k$ cells for the 3D domain. 
For the lower dimensional objects the number of cells should be chosen accordingly. 
The meshes may, if desirable, be obtained using the Gmsh mesh generator using the file \texttt{gmsh.geo} in the folder \texttt{regular/geometry} of the repository, where the mesh parameter 
\texttt{h} can vary in $[0.25, 0.12, 0.045]$ \si{\metre}. The time interval is divided in 100 equal time steps.
The results to be reported are:
\begin{enumerate} 
      
    \item The number of cells and degrees of freedom, according to the guidelines described in
    in Subsection \ref{subseq:data_reporting} point \ref{data:results}, saved in two files \texttt{results\_cond\_\$\{COND\}.csv};
      
	\item \label{case2_result2} for each mesh, the plot of the matrix hydraulic head $\head_3$ along the line from (0, 0, 0) to $(1,1,1)$;       
        
	\item \label{case2_result3} for the second level of mesh refinement, the average concentration on each matrix block during time. The average concentration is defined as $\int_{\Omega_{3,i}} c_3 / |\Omega_{3,i}|$\si{\meter\tothe{-3}}, with $i$ the region identifier.        
       
\end{enumerate}

In files called \texttt{dol\_cond\_\$\{COND\}\_refinement\_\$\{REFINEMENT\}.csv}
report the results from point \ref{case2_result2}. The variable \texttt{\$\{REFINEMENT\}} runs from 0 to 2 for increasing refinement, while \texttt{\$\{COND\}} should be 0 for the high-conductive case and 1 for low-conductive.
For point \ref{case2_result3}, the sub-domain identification indexes $region\_id$ are defined in Appendix \ref{appendix:regular}, and the file \texttt{color\_regions.vtu} shows a graphical representation of these regions.        The output is organised in the following way: two files named \texttt{dot\_cond\_\$\{COND\}.csv} having time as the first column, and the other columns the averaged concentration in each of the regions.

%%%%%%%%%%%%%%%%%%%%%%%%%%%%%%%%%%%%%%%%%%%%%%%%%%%%%%%%%%%%%%%%%%%%%%%%%%%%%%%%%%%%%%%%%%%%%%%%%%%%%%%%%%%%%%%%%%%%%%%%%%

\newpage
\subsection{Case 3: network with small features} \label{subsec:small_features}

\begin{figure}[!b]
\centering
\includegraphics[width=0.6\textwidth]{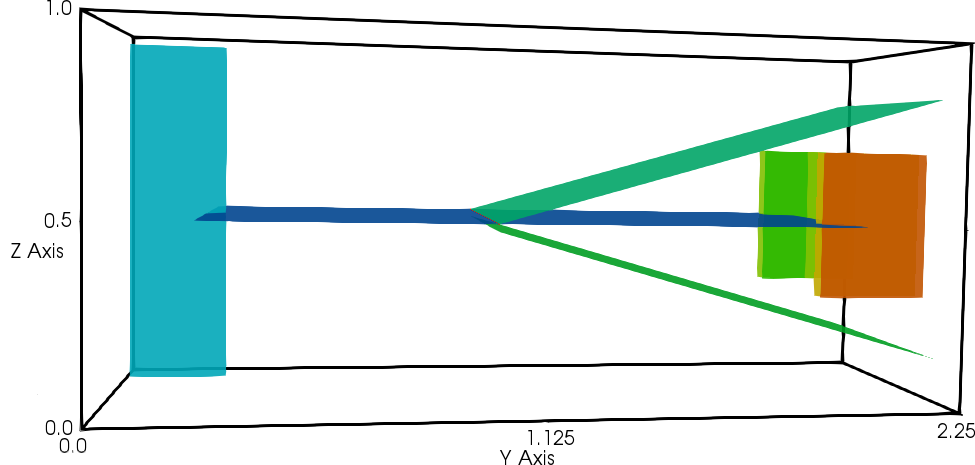}\hfill%
\includegraphics[width=0.3\textwidth]{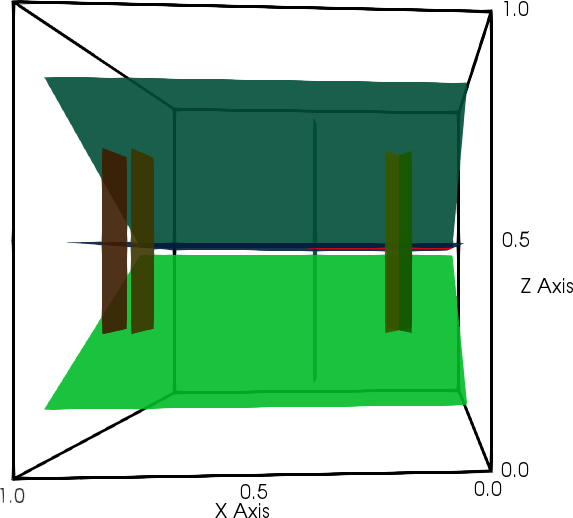}
\caption{Representation of the fractures and the outline of the domain
    	for the test case of Subsection \ref{subsec:small_features}.}
\label{fig:geom_small_features}
\end{figure}

\paragraph{Description}
This test is designed to probe accuracy in the presence of small geometric features, that may cause trouble for conforming meshing strategies.
The domain is a box $\Omega = (\SI{0}{\metre}, \SI{1}{\metre}) \times (\SI{0}{\metre}, \SI{2.25}{\metre}) \times (\SI{0}{\metre}, \SI{1}{\metre})$, containing eight fractures. These are defined by their vertexes, given in Appendix \ref{appendix:small_features}, and the geometry is depicted in Figure \ref{fig:geom_small_features}.
A Gmsh file \texttt{gmsh.geo} representing this geometry is uploaded to the Git repository in the folder \texttt{small\_features/geometry}.

We define in and outlet boundaries as follows:
\begin{gather*}
\partial\globalDomain_N = \partial\Omega \setminus (\partial\globalDomain_{in} \cup \partial\globalDomain_{out}) \qquad
\partial\globalDomain_{in} = (\SI{0}{\metre}, \SI{1}{\metre}) \times \{\SI{0}{\metre}\} \times (\SI{1/3}{\metre}, \SI{2/3}{\metre})\\
\partial\globalDomain_{out} = \partial\globalDomain_{out, 0} \cup \partial\globalDomain_{out, 1}\quad
\partial\globalDomain_{out, 0} = (\SI{0}{\metre}, \SI{1}{\metre}) \times \{\SI{2.25}{\metre}\} \times (\SI{0}{\metre}, \SI{1/3}{\metre}) \\
\partial\globalDomain_{out, 1} = (\SI{0}{\metre}, \SI{1}{\metre}) \times \{\SI{2.25}{\metre}\} \times (\SI{2/3}{\metre}, \SI{1}{\metre})
\end{gather*}
The boundary conditions for flow are: homogeneous Dirichlet conditions on $\partial\globalDomain_{out}$, uniform unit inflow on $\partial\globalDomain_{in}$, so that $\int_{\partial\globalDomain_{in}} \vecu_3 \cdot\normal dS = -1/3$ \si{\metre\cubed\per\second}, and homogeneous Neumann conditions on $\partial\globalDomain_N$. For the transport problem, we consider
homogeneous initial conditions and as boundary condition a unit concentration at $\partial\globalDomain_{in}$. The overall simulation time is \SI{1}{\second}.

\paragraph{Parameters}
\begin{center}
\begin{tabular}{|l|ll|}\hline
%\textbf{Parameter} & \multicolumn{2}{l|}{\textbf{Value}} \\\hline
Matrix hydraulic conductivity $K_3$& $\bm{I}$ & \si{\metre\per\second} \\
Fracture effective tangential hydraulic conductivity $K_2$ & \num{1e2}$\bm{I}$ & \si{\metre\squared\per\second} \\
Fracture effective normal hydraulic conductivity $\kappa_2$ & \num{2e6} & \si{\per\second} \\
Intersection effective tangential hydraulic conductivity $K_1$ & $1$ &  \si{\metre\cubed\per\second} \\
Intersection effective normal hydraulic conductivity $\kappa_1$ & \num{2e4} & \si{\metre\per\second} \\
Matrix porosity $\phi_3$ & \num{2e-1} & \\
Fracture porosity $\phi_2$ & \num{2e-1} & \\
Intersection effective porosity $\phi_1$ & \num{2e-1} &  \\
Fracture cross-sectional length $\epsilon_2$ & \num{1e-2} & \si{\metre} \\
Intersection cross-sectional area $\epsilon_1$ & \num{1e-4} & \si{\metre\squared} \\
\hline
\end{tabular}
\end{center}

\paragraph{Folder in Git repository} \texttt{small\_features}

\paragraph{Results}

The simulations should be done on 2 grids with approximately $30k$ and $150k$ 3D cells, labeled by \texttt{\$\{REFINEMENT\}} 0 and 1, respectively. The time step size is $\Delta t = \SI{1e-2}{s}$.
The results to be reported are:
\begin{enumerate} 
    \item \label{case3_result1} The number of cells and degrees of freedom according to the guidelines described in Subsection \ref{subseq:data_reporting} point \ref{data:results} in a file named \texttt{results.csv}. In addition,
 the total outflow over $\partial\globalDomain_{out, 0}$ and $\partial\globalDomain_{out, 1}$ as columns seven and eight; 
    \item \label{case3_result2} the average hydraulic head $\head_3$ at the inlet, $\partial\globalDomain_{in}$, reported as column nine of \texttt{results.csv}. The interpretation of hydraulic head at the inlet faces may differ between numerical methods; this should be seen as part of the discretization.
    \item \label{case3_result3} the hydraulic head $\head_3$ in the matrix along the line $\left(\SI{0.5}{\metre}, \SI{1.1}{\metre}, \SI{0}{\metre} \right)$-$\left(\SI{0.5}{\metre}, \SI{1.1}{\metre}, \SI{1}{\metre} \right)$;
    \item \label{case3_result4} the hydraulic head $\head_3$ in the matrix along $\left(\SI{0}{\metre}, \SI{2.15}{\metre}, \SI{0.5}{\metre} \right)$-$\left(\SI{1}{\metre}, \SI{2.15}{\metre}, \SI{0.5}{\metre} \right)$;
	\item \label{case3_result5} the mean tracer concentrations for each fracture and time step.    
\end{enumerate}
In files called \texttt{dol\_line\_\$\{LINE\}\_refinement\_\$\{REFINEMENT\}.csv}
report the results from point \ref{case3_result3} and \ref{case3_result4}, respectively for \texttt{\$\{LINE\} = 0} and \texttt{\$\{LINE\} = 1}. For point \ref{case3_result5},
report in a file \texttt{dot\_refinement\_\$\{REFINEMENT\}.csv}, with time in the first column, followed by the eight mean values ordered as in Appendix \ref{appendix:small_features}.

\newpage
\subsection{Case 4: a field case} \label{subsec:field_network}
%\cite{Matthai:2007:FEN}

\begin{figure}[!b]
\includegraphics[width=.55\textwidth]{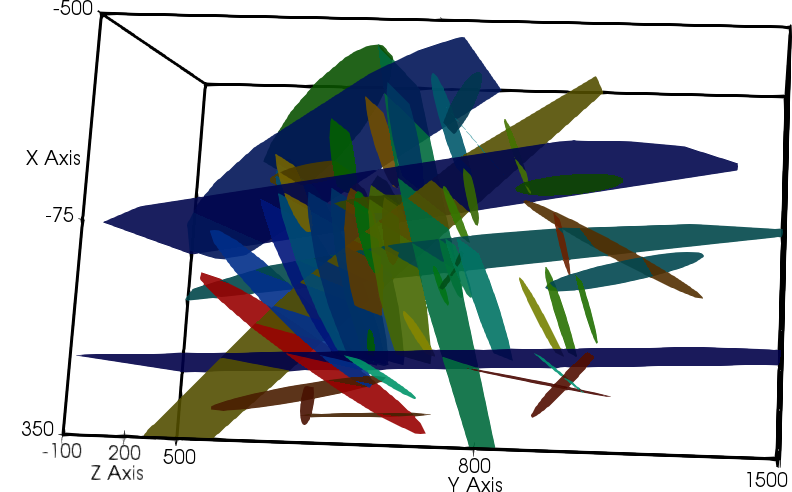}\hfill%
\includegraphics[width=.4\textwidth]{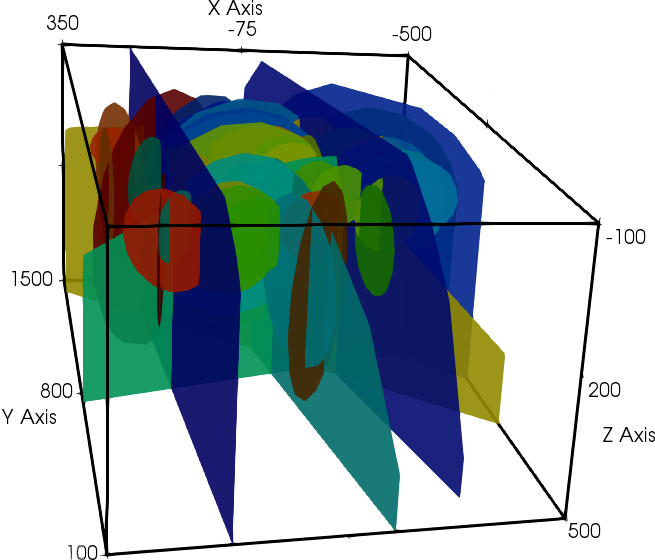}
\caption{Representation of the fractures and the outline of the domain
    	for the test case of Subsection \ref{subsec:field_network}.}
\label{fig:field_network}
\end{figure}

\paragraph{Description}
The fourth example has a geometry based on a post-processed %interpreted
outcrop from the island of Algerøyna, outside Bergen, Norway.
From the outcrop, 52 fractures are selected, extruded in the vertical direction and then cut by a bounding box. The resulting fracture geometry is depicted in Figure \ref{fig:field_network}.
The resulting network has 106 fracture intersections, and multiple fractures intersecting the domain boundary.

The simulation domain is the box $\Omega = (\SI{-500}{\metre}, \SI{350}{\metre}) \times (\SI{100}{\metre}, \SI{1500}{\metre}) \times (\SI{-100}{\metre}, \SI{500}{\metre})$. 
The polygons describing the fractures can be found in the file \texttt{fracture\_network.csv} in the folder \texttt{field}. A Gmsh file \texttt{gmsh.geo} describing the fracture geometry is provided in the subfolder \texttt{field/geometry}.
We define in and outlet boundaries as follows:
\begin{gather*}
\partial\globalDomain_N = \partial\Omega \setminus (\partial\globalDomain_{in} \cup \partial\globalDomain_{out}), \qquad
\partial\globalDomain_{in} = \partial\globalDomain_{in, 0} \cup \partial\globalDomain_{in, 1},\qquad
\partial\globalDomain_{out} = \partial\globalDomain_{out, 0} \cup \partial\globalDomain_{out, 1},\\
\partial\globalDomain_{in, 0} = (\SI{-500}{\metre}, \SI{-200}{\metre}) \times \{\SI{1500}{\metre}\} \times (\SI{300}{\metre}, \SI{500}{\metre}), \\
\partial\globalDomain_{in, 1} = \{\SI{-500}{\metre}\} \times (\SI{1200}{\metre}, \SI{1500}{\metre}) \times (\SI{300}{\metre}, \SI{500}{\metre}),\\
\partial\globalDomain_{out, 0} = \{\SI{-500}{\metre}\} \times (\SI{100}{\metre}, \SI{400}{\metre}) \times (\SI{-100}{\metre}, \SI{100}{\metre}),\\
\partial\globalDomain_{out, 1} = \{\SI{350}{\metre}\} \times (\SI{100}{\metre}, \SI{400}{\metre}) \times (\SI{-100}{\metre}, \SI{100}{\metre}).
\end{gather*}
The boundary conditions for flow are: homogeneous Dirichlet conditions on $\partial\globalDomain_{out}$, uniform unit inflow on $\partial\globalDomain_{in}$, so that $\int_{\partial\globalDomain_{in}} \vecu_3 \cdot\normal dS = \SI{-1.2e5}{\metre\cubed\per\second}$, and homogeneous Neumann conditions on $\partial\globalDomain_N$. For the transport problem, we consider
homogeneous initial conditions, with a unit concentration as boundary condition at $\partial\globalDomain_{in}$. The overall simulation time is \SI{5e3}{\second}.

\paragraph{Parameters}
\begin{center}
\begin{tabular}{|l|ll|}\hline
%\textbf{Parameter} & \multicolumn{2}{l|}{\textbf{Value}} \\\hline
Matrix hydraulic conductivity $K_3$& $\bm{I}$ & \si{\metre\per\second} \\
Fracture effective tangential hydraulic conductivity $K_2$ & \num{1e2}$\bm{I}$ & \si{\metre\squared\per\second} \\
Fracture effective normal hydraulic conductivity $\kappa_2$ & \num{2e6} & \si{\per\second} \\
Intersection effective tangential hydraulic conductivity $K_1$ & $1$ &  \si{\metre\cubed\per\second} \\
Intersection effective normal hydraulic conductivity $\kappa_1$ & \num{2e4} & \si{\metre\per\second} \\
Matrix porosity $\phi_3$ & \num{2e-1} & \\
Fracture porosity $\phi_2$ & \num{2e-1} & \\
Intersection porosity $\phi_1$ & \num{2e-1} & \\
Fracture cross-sectional length $\epsilon_2$ & \num{1e-2} & \si{\metre} \\
Intersection cross-sectional area $\epsilon_1$ & \num{1e-4} & \si{\metre\squared} \\
\hline
\end{tabular}
\end{center}

\paragraph{Folder in Git repository} \texttt{field}

\paragraph{Results}
We consider a time step size of $\Delta t = \SI{50}{\second}$. The results to be reported are:
\begin{enumerate} 
  \item \label{case4_result1} The number of cells and degrees of freedom according to the guidelines described in Subsection \ref{subseq:data_reporting} point \ref{data:results} in a file named \texttt{results.csv}. In addition,
 the total outflow over $\partial\globalDomain_{out, 0}$ and $\partial\globalDomain_{out, 1}$ as columns seven and eight;
 
    \item \label{case4_result2} the average hydraulic head $\head_3$ at the inlet, $\partial\globalDomain_{in}$, reported as the seventh column of \texttt{results.csv};
    
    \item\label{case4_result3} the hydraulic head $\head_3$ in the matrix along $\left(\SI{350}{\metre}, \SI{100}{\metre}, \SI{-100}{\metre}\right)$-$\left(\SI{-500}{\metre}, \SI{1500}{\metre}, \SI{500}{\metre} \right)$;
    
    \item\label{case4_result4} the hydraulic head $\head_3$ in the matrix along $\left(\SI{-500}{\metre}, \SI{100}{\metre}, \SI{-100}{\metre}\right)$-$\left(\SI{350}{\metre}, \SI{1500}{\metre}, \SI{500}{\metre} \right)$;
        
	\item\label{case4_result5} the mean tracer concentrations throughout the simulation time in all fractures.
\end{enumerate}
In files called \texttt{dol\_line\_\$\{LINE\}.csv}
report the results from point \ref{case4_result3} and \ref{case4_result4}, respectively for \texttt{\$\{LINE\} = 0} and \texttt{\$\{LINE\} = 1}. For point \ref{case4_result5},
report in a file \texttt{dot.csv}, with time in the first column, followed by the 52 mean concentration values.

\newpage
\section{How to participate}
\label{sec:participation}
We welcome all computational modelers to participate in the proposed
benchmark study. Three different stages of participation can be identified:
1. registration, 2. simulation and synchronization, 3. publication.
In the following, we will describe each stage in more detail.

\subsection{Registration}
If an interested researcher would like to participate in the benchmark study,
he/she should fill in the Google form
provided at \url{https://goo.gl/forms/wCaPLipPI1d5Qz1j1}. One form should be submitted for
each discretization method. The form asks for personal information from the participant(s)
as well as for a brief description of the discretization method.
After the application has been approved, the researcher will be asked to sign a
participation agreement that states that the participant should not publish the results
of other people before we do collectively.
Once that is signed, the researcher will be granted developer access to the Git repository
enabling him/her to download the required input data for the benchmark cases,
to view the simulation results of the other participants
as well as to upload his/her own results.

\subsection{Simulation and synchronization}
In this stage, the participant is expected to perform the computations and convert
his/her results to the
requested format. In particular, he/she should try the proposed scripts to evaluate the data
as indicated in each case description.
The deadline for uploading the results to the Git repository is 31 January 2019.
The participant should follow the GitLab feature branch workflow that consists of
contributing additions via a GitLab merge request of a corresponding branch,
analogous to the
GitHub standard fork and pull request workflow, see the details below in Subsection \ref{subseq:data_reporting}.
The uploaded results will be visible to all participants. This will initiate a synchronization
phase in which we will discuss and potentially recompute/adjust the results.
In particular, a minisymposium at the SIAM GS 19 is planned for presenting the benchmark cases, the participants and the participating schemes as well as for discussing the results. The synchronization phase is supposed to be finished by uploading the final results before 15 May 2019.

\subsection{Publication}
As a final outcome of this benchmark study, we aim for a joint publication in a peer-reviewed
top-ranked journal. We would expect from each participant that he/she helps us in writing
and evaluating the manuscript before initial submission as well as during the reviewing process.
We would like to stress that co-authorship can only be granted for at most two persons
per participating discretization method.

\subsection{Details on reporting results}
\label{subseq:data_reporting}
All data requested will be collected in the Git repository \url{{https://git.iws.uni-stuttgart.de/benchmarks/fracture-flow-3d.git}}. To ensure a unified and streamlined work flow, please observe the following structure and formatting of the uploaded files.

The repository contains one folder for each of the four test cases. In the \texttt{results} folder of these, each method will be assigned a folder with the format
\begin{center}
	\texttt{\$ACRONYM\_INSTITUTION/\$ACRONYM\_NUMERICAL\_SCHEME/}. 
\end{center}
No white spaces are allowed, please be consistent with upper and lower case.
As an example, results for the first test case obtained using the University of Bergen implementation of the Mixed Virtual Element method are stored in \texttt{\seqsplit{single/results/UiB/MVEM}}.
Similarly, the string \texttt{\$ACRONYM\_INSTITUTION-\$ACRONYM\_NUMERICAL\_SCHEME}  (e.g. \texttt{UiB-MVEM}) constitutes an \texttt{ID} to be used in the handling and labeling of the results.

In the methods' folder, three categories of files should be placed, for which we here specify the format:
\begin{enumerate}
\item \label{data:pol} \textbf{Plot over line:} A csv file whose first column contains the arc length and the second column contains the corresponding values for either $c$ or $\head$, sampled at 2k equidistant points along the specified line. The delimiter is a comma and there is no header for the file, an example is
\begin{verbatim}
0.0,3.36948657742311
0.5,5.38490345323433
1.0,8.34820934803293
\end{verbatim}
Generally, the file is called \texttt{dol} (data over line). If several refinement levels are reported for the same test case, the files will be distinguished by adding an appendix to the file name, see individual test case descriptions. The plot over line can be easily done with the ``Plot Over Line'' filter of ParaView, but a rearrangement of the output may be required.

\item \label{data:pot} \textbf{Plot versus time:} A csv file named \texttt{dot} (data over time) whose first column contains the time and the subsequent columns contain the values specified in each test case. The same specifications apply as in the previous case.

\item \label{data:results} \textbf{Cell numbers and other metrics:} A file named \texttt{results.csv} where the first four columns contain the total number of cells in $\globalDomain_d, \ d = 0, \dots, 3$. The next two columns contain the number of degrees of freedom \texttt{dof} and the number of non-zero elements in the matrix \texttt{nnz}, respectively. 
These are followed by other metrics in the order specified in the test case descriptions. If results are reported for several refinement levels, the coarsest refinement level corresponds to the first row etc. Again, the delimiter is a comma and there is no header for the file.
\end{enumerate}

Scripts for generating the plots are to be found in the \texttt{scripts} folders in each of the test case folders. The scripts are named \textit{pol.py} and \textit{pot.py} and contain comments on how to add new data. Once new results have been generated, first update the Git repository. Add corresponding lines to the plot scripts and upload the scripts and commit the corresponding data files to a newly created feature branch. After pushing to the remote server, you can create a corresponding merge request that will be reviewed and eventually approved by the repository maintainers. Please verify that the scripts are working before pushing to the remote Git repository. More details can be found at \url{https://docs.gitlab.com/ee/workflow/workflow.html} and \url{https://docs.gitlab.com/ee/gitlab-basics/add-merge-request.html}.

\section{Appendix}
\subsection{Additional information for the test case of Subsection \ref{subsec:regular}} \label{appendix:regular}
We report the coordinates for the four corners of the fractures
\begin{center}
	\begin{tabular}{|c|c c c| c c c| c c c| c c c|}
		\hline
        $ID$ & $x_0$ & $y_0$ & $z_0$ & $x_1$ & $y_1$ & $z_1$ & $x_2$ & $y_2$ & $z_2$ & $x_3$ & $y_3$ & $z_3$\\ \hline
		0 & 0.5 & 0 & 0 & 0.5 & 1 & 0 & 0.5 & 1 & 1 & 0.5 & 0 & 1 \\ \hline                  
		1 & 0 & 0.5 & 0 & 1 & 0.5 & 0 & 1 & 0.5 & 1 & 0 & 0.5 & 1 \\ \hline                  
		2 & 0 & 0 & 0.5 & 1 & 0 & 0.5 & 1 & 1 & 0.5 & 0 & 1 & 0.5 \\ \hline                  
		3 & 0.75 & 0.5 & 0.5 & 0.75 & 1.0 & 0.5 & 0.75 & 1.0 & 1.0 & 0.75 & 0.5 & 1.0 \\ \hline
		4 & 0.5 & 0.5 & 0.75 & 1.0 & 0.5 & 0.75 & 1.0 & 1.0 & 0.75 & 0.5 & 1.0 & 0.75 \\ \hline
		5 & 0.5 & 0.75 & 0.5 & 1.0 & 0.75 & 0.5 & 1.0 & 0.75 & 1.0 & 0.5 & 0.75 & 1.0 \\ \hline
		6 & 0.50 & 0.625 & 0.50 & 0.75 & 0.625 & 0.50 & 0.75 & 0.625 & 0.75 & 0.50 & 0.625 & 0.75 \\ \hline
		7 & 0.625 & 0.50 & 0.50 & 0.625 & 0.75 & 0.50 & 0.625 & 0.75 & 0.75 & 0.625 & 0.50 & 0.75 \\ \hline
		8 & 0.50 & 0.50 & 0.625 & 0.75 & 0.50 & 0.625 & 0.75 & 0.75 & 0.625 & 0.50 & 0.75 & 0.625 \\ \hline
	\end{tabular}
\end{center}

\noindent We report the $region\_id$ 
       	\begin{gather*}
        	region\_id(x<0.5 \cap y<0.5 \cap z<0.5) = 0\\                                             
			region\_id(x>0.5 \cap y<0.5 \cap z<0.5) = 1\\                                             
			region\_id(x<0.5 \cap y>0.5 \cap z<0.5) = 2\\                                             
			region\_id(x>0.5 \cap y>0.5 \cap z<0.5) = 3\\                                             
			region\_id(x<0.5 \cap y<0.5 \cap z>0.5) = 4\\                                             
			region\_id(x>0.5 \cap y<0.5 \cap z>0.5) = 5\\                                             
			region\_id(x<0.5 \cap y>0.5 \cap z>0.5) = 6 \\           
			region\_id(x>0.75 \cap y>0.75 \cap z>0.75) = 7 \\                                         
			region\_id(x>0.75 \cap y>0.5 \cap y<0.75 \cap z>0.75) = 8 \\                             
			region\_id(x>0.5 \cap x<0.75 \cap y>0.75 \cap z>0.75) = 9 \\                             
			region\_id(x>0.5 \cap x<0.75 \cap y>0.5 \cap y<0.75 \cap z>0.75) = 10  \\               
			region\_id(x>0.75 \cap y>0.75 \cap z>0.5 \cap z<0.75) = 11 \\                            
			region\_id(x>0.75 \cap y>0.5 \cap y<0.75 \cap z>0.5 \cap z<0.75) = 12  \\               
			region\_id(x>0.5 \cap x<0.75 \cap y>0.75 \cap z>0.5 \cap z<0.75) = 13 \\          
			region\_id(x>0.5 \cap x<0.625 \cap y>0.5 \cap y<0.625 \cap z>0.5 \cap z<0.625) = 14  \\
			region\_id(x>0.625 \cap x<0.75 \cap y>0.5 \cap y<0.625 \cap z>0.5 \cap z<0.625) = 15 \\
			region\_id(x>0.5 \cap x<0.625 \cap y>0.625 \cap y<0.75 \cap z>0.5 \cap z<0.625) = 16 \\
			region\_id(x>0.625 \cap x<0.75 \cap y>0.625 \cap y<0.75 \cap z>0.5 \cap z<0.625) = 17\\
			region\_id(x>0.5 \cap x<0.625 \cap y>0.5 \cap y<0.625 \cap z>0.625 \cap z<0.75) = 18 \\
			region\_id(x>0.625 \cap x<0.75 \cap y>0.5 \cap y<0.625 \cap z>0.625 \cap z<0.75) = 19\\
			region\_id(x>0.5 \cap x<0.625 \cap y>0.625 \cap y<0.75 \cap z>0.625 \cap z<0.75) = 20\\
			region\_id(x>0.625 \cap x<0.75 \cap y>0.625 \cap y<0.75 \cap z>0.625 \cap z<0.75) = 21
        \end{gather*}

\subsection{Additional information for the test case of Subsection \ref{subsec:small_features}} \label{appendix:small_features}
We report the coordinates for the four corners of the fractures
\begin{center}
	\begin{tabular}{|c|c c c| c c c| c c c| c c c|}
		\hline
        $ID$ & $x_0$ & $y_0$ & $z_0$ & $x_1$ & $y_1$ & $z_1$ & $x_2$ & $y_2$ & $z_2$ & $x_3$ & $y_3$ & $z_3$\\ \hline
		0 & 0.05 & 0.25 & 0.5 & 0.95 & 0.25 & 0.5 & 0.95 & 2 & 0.5 & 0.05 & 2 & 0.5 \\ \hline                  
		1 & 0.5 & 0.05 & 0.05 & 0.5 & 0.05 & 0.95 & 0.5 & 0.3 & 0.95 & 0.5 & 0.3 & 0.05 \\ \hline                  
		2 & 0.05 & 1 & 0.5 & 0.95 & 1 & 0.5 & 0.95 & 2.2 & 0.85 & 0.05 & 2.2 & 0.85 \\ \hline                  
		3 & 0.05 & 1 & 0.48 & 0.95 & 1.0 & 0.48 & 0.95 & 2.2 & 0.14 & 0.05 & 2.2 & .14 \\ \hline
		4 & 0.17 & 1.9 & 0.7 & 0.17 & 1.9 & 0.3 & 0.23 & 2.2 & 0.3 & 0.23 & 2.2 & 0.7 \\ \hline
        5 & 0.23 & 1.9 & 0.7 & 0.23 & 1.9 & 0.3 & 0.17 & 2.2 & 0.3 & 0.17 & 2.2 & 0.7 \\ \hline
		6 & 0.77 & 1.9 & 0.7 & 0.77 & 1.9 & 0.3 & 0.77 & 2.2 & 0.3 & 0.77 & 2.2 & 0.7 \\ \hline
        7 & 0.83 & 1.9 & 0.7 & 0.83 & 1.9 & 0.3 & 0.83 & 2.2 & 0.3 & 0.83 & 2.2 & 0.7 \\ \hline
	\end{tabular}
\end{center}

\bibliography{literature,biblio}

%\section*{Acknowledgements (not compulsory)}

%Acknowledgements should be brief, and should not include thanks to anonymous referees and editors, or effusive comments. Grant or contribution numbers may be acknowledged.

%\section*{Author contributions statement}

%Must include all authors, identified by initials, for example:
%A.A. conceived the experiment(s),  A.A. and B.A. conducted the experiment(s), C.A. and D.A. analysed the results.  All authors reviewed the manuscript. 

%\section*{Additional information}

%To include, in this order: \textbf{Accession codes} (where applicable); \textbf{Competing financial interests} (mandatory statement). 

%The corresponding author is responsible for submitting a \href{http://www.nature.com/srep/policies/index.html#competing}{competing financial interests statement} on behalf of all authors of the paper. This statement must be included in the submitted article file.
\end{document}